\documentclass{amsart}
\usepackage{amssymb} 
\newtheorem{thm}{Theorem}[section]
\newtheorem{cor}[thm]{Corollary}

\newtheorem{prop}[thm]{Proposition}

\newtheorem{rem}[thm]{Remark}

\theoremstyle{question}
\newtheorem{qu}[thm]{Question}

\numberwithin{equation}{section} 

\begin{document}
\author[Abdollahi and Khosravi]{A. Abdollahi \& H. Khosravi}
\title[Right $n$-Engel elements]
{When right $n$-Engel elements of a group form a subgroup?}
\address{Department of Mathematics, University of Isfahan, Isfahan 81746-73441, Iran}%
\email{alireza\_abdollahi@yahoo.com}
\email{hassan\_khosravy@yahoo.com}
\thanks{}
\begin{abstract}
Let $R_n(G)$ denotes the set of all  right $n$-Engel elements of a
group $G$. We show  that in any group $G$ whose $5$th term of lower central series has no  element of
order $2$, $R_3(G)$ is a subgroup. Furthermore we prove that
$R_4(G)$ is a subgroup for locally nilpotent groups $G$ without
elements of orders $2$, $3$ or $5$; and in this case the normal
closure $\langle x \rangle^G$ is nilpotent of class at most $7$
for each $x\in R_4(G)$. Using a group constructed by Newman and
Nickel we also show that,   for each $n\geq 5$, there exists a
nilpotent group of class $n+2$  containing a right $n$-Engel
element $x$ and an element $a\in G$ such that both $[x^{-1},_n a]$ and $[x^{k},_n a]$ are of infinite order for all integers $k\geq 2$. We finish the paper by proving that at least one of the following  happens:
(1) \; There is an infinite finitely generated $k$-Engel group of exponent $n$ for some positive integer $k$ and some $2$-power number $n$.
(2) \; There is a  group  generated by finitely  many bounded left Engel elements which is not  an Engel group.
\end{abstract}
\maketitle
\section{\bf Introduction and Results}\label{sec1}
 Let $G$ be any group and
$n$  a non-negative integer. For any two
 elements $a$ and $b$ of $G$, we define
inductively $[a,_n b]$ the $n$-Engel commutator of the pair
$(a,b)$, as follows:
$$[a,_0 b]:=a,~ [a,b]=[a,_1 b]:=a^{-1}b^{-1}ab \mbox{ and } [a,_n
b]=[[a,_{n-1} b],b]\mbox{ for all }n>0.$$  An element $x$ of $G$
is called right (left, resp.) $n$-Engel if $[x,_ng]=1$ ($[g,_n
x]=1$, resp.) for all $g\in G$. We denote by $R_n(G)$ ($L_n(G)$,
resp.), the set of all right (left, resp.) $n$-Engel elements of
$G$.  A group $G$ is called $n$-Engel if $G=L_n(G)$ or
equivalently $G=R_n(G)$. It is clear that $R_0(G)=1$,
$R_1(G)=Z(G)$ the center of $G$, and  W.P. Kappe
\cite{kappew} (implicitly) proved $R_2(G)$ is a characteristic
subgroup of $G$. L.C. Kappe and Ratchford \cite{kappe2} have shown that $R_n(G)$ is a subgroup of $G$ whenever
 $G$ is a metabelian group, or $G$ is a
center-by-metabelian group such that $[\gamma_k(G),\gamma_j(G)]=1$
for some $k,j\geq 2$ with $k+j-2\leq n$ and $n\geq 3$. Macdonald \cite{macd} has shown
that the inverse or square of a right 3-Engel element need not be
right 3-Engel. Nickel \cite{nick2} generalized Macdonald's result to all $n
\geq 3$. In fact he constructed a group with a right $n$-Engel
element $a$  neither $a^{-1}$ nor $a^2$ is a right $n$-Engel
element. The construction of Nickel's example was guided by
computer experiments and  arguments  based on commutator
calculus. Although Macdonald's example shows that $R_3(G)$ is not
in general a subgroup of $G$, Heineken \cite{hein} has already shown that
if $A$ is the subset of a group $G$ consisting of all elements $a$
such that $a^{\pm 1}\in R_3(G)$, then $A$ is a subgroup if either
$G$ has no element of order $2$ or $A$ consists only of elements
having finite odd order. Newell \cite{newell} proved that the
normal closure of every right $3$-Engel element is nilpotent of
class at most $3$. In Section 2 we prove that if $G$ is a
$2'$-group, then $R_3(G)$ is a subgroup of $G$. Nickel's example
shows that the set of right 4-Engel elements is not a subgroup in
general (see also  first Example in Section 4 of \cite{ab1}). In
Section 3, we prove that if $G$ is a locally nilpotent
$\{2,3,5\}'$-group, then $R_4(G)$ is a subgroup of $G$.

Traustason \cite{traus} proved that any locally nilpotent 4-Engel
group $H$ is Fitting of degree at most $4$. This means that the
normal closure of every element of $H$ is nilpotent of class at
most 4. More precisely he proved that if $H$ has no element of
order $2$ or $5$, then $H$ has Fitting degree at most $3$. Now by a result of Havas and Vaughan-Lee \cite{havas}, one knows  any 4-Engel group is  locally nilpotent and so Traustason's  result is true for all $4$-Engel groups. In Section 3, by another result of Traustason \cite{traus2} we show that the normal closure of every right $4$-Engel element in a locally nilpotent
$\{2,3,5\}'$-group, is nilpotent of class at most $7$.

Newman and Nickel \cite{newman} have shown that for every $n\geq
5$ there exists a nilpotent group $G$ of class $n+2$ containing a
right $n$-Engel element $a$ and an element $b$ such that $[b,_na]$
has infinite order. As we mentioned  above, Nickel
\cite{nick2} has shown that for every $n\geq 3$ there exists a
nilpotent group of class $n+2$ having a right $n$-Engel element
$a$ and an element $b$ such that $[a^{-1},_n b]=[a^2,_n b]\neq 1$. We have checked that  the latter element in Nickel's example
 is of finite order whenever $n\in\{5,6,7,8\}$. In Section \ref{se4}, using
the group constructed by Newman and Nickel we show
that there exists a nilpotent group $G$ of class $n+2$ such that
$x\in R_n(G)$ and both $[x^{-1},_n a]$ and $[x^{k},_na]$ have infinite order for every integer
$k\geq 2$.\\

In \cite{ab1} the following question has been proposed:
\begin{qu}\label{q1}
 Let $n$ be a positive integer. Is there a set of prime
 numbers $\pi_n$ depending only on $n$ such that the set
 of right $n$-Engel elements in any nilpotent or finite
 $\pi'_n$-group forms a subgroup?
\end{qu}
In Section 4  we negatively answers
Question \ref{q1}.\\

  As far as we know there is no published example of a group  whose set of (bounded) right Engel elements do not form a subgroup.  But for the set  of bounded left Engel elements there are some evidences supporting this idea that the subgroupness of bounded left Engel elements of a an arbitrary group  should be false.  We finish the paper by proving that   at least one of the following  happens:
 \begin{enumerate}
 \item There is an  infinite finitely generated $k$-Engel group of exponent $n$ for some positive integer $k$ and some $2$-power number $n$.
 \item There is a  group  generated by finitely  many bounded left Engel elements which is not  an Engel group, where by an Engel group we mean a group in which for every two elements $x$ and $y$, there exists an integer $k=k(x,y)\geq 0$ such that $[x,_k y]=1$.
 \end{enumerate}

 Throughout the paper we have frequently use {\sf GAP nq}
package of Werner Nickel.
 All given  timings  were obtained on an Intel
Pentium 4-1.70GHz processor with 512 MB running Red Hat Enterprise Linux 5.\\
\section{\bf Right 3-Engel elements}\label{se2}
Throughout for any positive integer $k$ and any group $H$,  $\gamma_k(H)$ denotes the $k$th term of the lower central series of $H$.
The main result of  this section implies that  $R_3(G)$ is a subgroup of $G$ whenever $G$ is a $2'$-group.
Newell \cite{newell} proved that
\begin{thm}\label{th1}
Let $G=\langle a,b,c\rangle$ be a group such that $a,b\in R_3(G)$.
Then
\begin{enumerate}
    \item $\langle a,c\rangle$ is nilpotent of class at most $5$ and
    $\gamma_5(\langle a,c\rangle)$ has exponent $2$.
    \item $G$ is nilpotent of class at most $6$.
    \item $\frac{\gamma_5(G)}{\gamma_6(G)}$ has exponent $10$.
    Furthermore $[a,c,b,c,c]^2\in \gamma_6(G)$.
    \item $\gamma_6(G)$ has exponent $2$.
\end{enumerate}
\end{thm}
\begin{thm}\label{th0}
Let $G$ be a group such that $\gamma_5(G)$ has no element of order $2$. Then $R_3(G)$ is a subgroup of $G$.
\end{thm}
\begin{proof}
Let $a,b\in R_3(G)$ and $c\in G$. We first show that $a^{-1}\in R_3(G)$. We have
\begin{eqnarray}
[a^{-1},c,c,c]&=&[[[a,c,a^{-1}]^{-1}[a,c]^{-1},c],c]\nonumber\\
&=&[[a,c,a,a^{-1},c][a,c,a,c],c][[a,c,c,[a,c]^{-1}][a,c,c]^{-1},c]\nonumber\\
&=&[a,c,a,c,c][a,c,c,c]^{[a,c,c]}\nonumber\\
&=&[a,c,a,c,c]\nonumber
\end{eqnarray}
Therefore by Theorem \ref{th1} (2), $a^{-1}\in R_3(G)$. On the
other hand
\begin{eqnarray}
[ab,c,c,c]&=&[[a,c][a,c,b][b,c],c,c]\nonumber\\
&=&[[a,c,c][a,c,c,[a,c,b]][[a,c,c],[b,c]][a,c,b,c]
[b,c,b,c,[b,c]][b,c,c],c]\nonumber\\
&=&[a,c,c,[b,c],c][a,c,b,c,c].\nonumber
\end{eqnarray}
Now by Theorem \ref{th1} $[a,c,c,[b,c],c],~[a,c,b,c,c]^2 \in
\gamma_6(G)$ and thus $ab\in R_3(G)$.
\end{proof}
Now we give a proof of Theorem \ref{th0} by using {\sf GAP nq}
package of Werner Nickel.

\noindent{\bf Second Proof of Theorem \ref{th0}.} By Theorem
\ref{th1}, we know that $\langle x,y,z\rangle$ is nilpotent if
$x,y\in R_3(G)$ and $z\in G$. We now construct the largest
nilpotent group $H=\langle a,b,c \rangle$ such that $a,b\in
R_3(H)$ and $c\in H$, by {\sf nq} package.
\begin{verbatim}
 LoadPackage("nq"); #nq package of Werner Nickel #
 F:=FreeGroup(4);a1:=F.1; b1:=F.2; c1:=F.3; x:=F.4;
 L:=F/[LeftNormedComm([a1,x,x,x]),LeftNormedComm([b1,x,x,x])];
 H:=NilpotentQuotient(L,[x]);
 a:=H.1; b:=H.2; c:=H.3;  d:=LeftNormedComm([a^{-1},c,c,c]);
 e:=LeftNormedComm([a*b,c,c,c]);  Order(d); Order(e);
 C:=LowerCentralSeries(H);  d in C[5]; e in C[5];
\end{verbatim}
Then if we consider the elements $d=[a^{-1},c,c,c]$ and
$e=[ab,c,c,c]$ of $H$, we can see by above command in {\sf GAP} that $d$ and $e$ are elements of $\gamma_5(H)$ and have orders $2$ and
$4$, respectively. So, in the group $G$, we have   $d=e=1$.
This completes the proof. $\hfill \Box$\\

Note that, the second proof of Theorem \ref{th0} also shows the necessity of assuming  that $\gamma_5(G)$ has no element of order $2$.
\section{\bf Right 4-Engel elements}\label{se3}
Our main result in this section is to prove the following.
\begin{thm}\label{th2}
Let $G$ be a  $\{2,3,5\}'$-group such that $\langle a,b,x\rangle$ is nilpotent for all $a,b\in R_4(G)$ and any $x\in G$. Then $R_4(G)$
is a subgroup of $G$.
\end{thm}
\begin{proof}
Consider the `freest' group, denoted by $U$, generated by two elements $u$,$v$
with $u$ a right 4-Engel element. We mean this by the group $U$
given by the presentation
$$\langle u,v \;|\; [u,_4 x]=1 \;\;\text{for all words}\;\; x \in F_2\rangle,$$
where $F_2$ is the free group generated by $u$ and $v$.
We do not know whether $U$ is nilpotent or not.
Using the {\sf nq} package shows that the group
$U$ has a largest nilpotent quotient $M$ with
class $8$.
By the following code, the
group $M$ generated by a right $4$-Engel element $a$ and an
arbitrary element $c$ is constructed.
We then see that the element $[a^{-1},c,c,c,c]$ of $M$
is of order  $375=3\times 5^3$.  Therefore the inverse of a right
$4$-Engel element of $G$ is again a right $4$-Engel element. The
following code in {\sf GAP} gives a proof of the latter claim. The computation
was completed in about 24 seconds.
\begin{verbatim}
 F:=FreeGroup(3); a1:=F.1; b1:=F.2; x:=F.3;
 U:=F/[LeftNormedComm([a1,x,x,x,x])];
 M:=NilpotentQuotient(U,[x]);
 a:=M.1; c:=M.2;
 h:=LeftNormedComm([a^-1,c,c,c,c]);
 Order(h);
\end{verbatim}
We now show that the product of every two
right 4-Engel elements in $G$
is a right 4-Engel element. Let $a,b\in R_4(G)$ and $c\in G$. Then
we claim that
$$H=\langle a,b,c\rangle \;\; \text{is nilpotent of class at most}\; 7. \;\;\;(*)$$ By induction on the nilpotency class of $H$, we may assume that $H$ is
nilpotent of class at most 8. Now we construct the largest nilpotent group
$K=\langle a_1,b_1,c_1\rangle$ of class 8 such that $a_1,b_1\in R_4(K)$.
\begin{verbatim}
 F:=FreeGroup(4);A:=F.1; B:=F.2; C:=F.3; x:=F.4;
 W:=F/[LeftNormedComm([A,x,x,x,x]),LeftNormedComm([B,x,x,x,x])];
 K:=NilpotentQuotient(W,[x],8);
 LowerCentralSeries(K);
\end{verbatim}
The computation took about 22.7 hours. We see that $\gamma_8(K)$
has exponent $60$. Therefore, as $H$ is a  $\{2,3,5\}'$-group, we have
$\gamma_8(H)=1$ and this completes the proof of our  claim $(*)$. \\
Therefore we have proved that any nilpotent group without elements of orders $2$, $3$ or $5$ which is generated by three elements two of which are right $4$-Engel, is nilpotent of class at most $7$.\\
Now we construct, by the {\sf nq} package, the largest nilpotent group $S$ of class $7$ generated by two right $4$-Engel elements $s,t$ and an arbitrary element $g$.  Then one can find  by {\sf GAP} that the order of   $[st,g,g,g,g]$ in $S$  is
  300. Since $H$ is a  quotient of $S$, we have that $[ab,c,c,c,c]$ is of order dividing $300$ and so it is trivial, since $H$ is a $\{2,3,5\}'$-group.
 This completes the proof.
\end{proof}
\begin{cor}\label{co1}
Let $G$ be a  $\{2,3,5\}'$-group such that $\langle a,b,x\rangle$ is nilpotent for all $a,b\in R_4(G)$ and for any $x\in G$. Then $R_4(G)$ is a nilpotent group of class at most $7$. In particular,  the normal closure of every right $4$-Engel element of group $G$ is nilpotent
of class at most $7$.
\end{cor}
\begin{proof}
By  Theorem \ref{th2}, $R_4(G)$ is a subgroup of $G$ and so it
is a 4-Engel group. In \cite{traus2} it is shown that every
locally nilpotent 4-Engel $\{2,3,5\}'$-group is nilpotent of class at most 7.
Therefore $R_4(G)$ is nilpotent of class at most 7. Since $R_4(G)$ is a normal set, the second part follows easily.
\end{proof}
Therefore, to prove that the normal closure of any right $4$-Engel element of a $\{2,3,5\}'$-group $G$ is nilpotent, it is enough to show that
$\langle a,b,x\rangle$ is nilpotent for all $a,b\in R_4(G)$ and for any $x\in G$. It may be  surprising that  Newell \cite{newell} has had a similar obstacle to prove that the normal closure of a right $3$-Engel element is nilpotent in any group.
\begin{cor}
In any $\{2,3,5\}'$-group, the normal closure of any right $4$-Engel element is nilpotent if and only if every $3$-generator subgroup in which two of the generators can be chosen to be  right $4$-Engel, is nilpotent.
\end{cor}
\begin{proof}
By Corollary \ref{co1}, it is enough to show that  a $\{2,3,5\}'$-group $H=\langle a,b,x\rangle$ is nilpotent whenever $a,b\in R_4(H)$, $x\in H$ and both $\langle a\rangle^H$ and $\langle b\rangle ^H$ are nilpotent. Consider the subgroup $K=\langle a\rangle ^H\langle b\rangle^H$ which is nilpotent by Fitting's theorem. Now we prove that $K$ is finitely generated. We have $K=\langle a,b\rangle^{\langle x\rangle}$ and since $a$ and $b$ are both right $4$-Engel, it is well-known that
$$\langle a\rangle^{\langle x\rangle}=\langle a,a^x,a^{x^2},a^{x^3}\rangle \;\;\text{and}\;\; \langle b\rangle^{\langle x\rangle}=\langle b,b^x,b^{x^2},b^{x^3}\rangle,$$
and so $$K=\langle a,a^x,a^{x^2},a^{x^3},b,b^x,b^{x^2},b^{x^3} \rangle.$$
It follows that $H$  satisfies maximal condition on its subgroups as it is (finitely generated nilpotent)-by-cyclic. Now by a famous result of Baer \cite{Baer} we have that $a$ and $b$ lie in the $(m+1)$th term $\zeta_m(H)$ of the upper central series of $H$ for some positive integer $m$. Hence $H/\zeta_m(H)$ is cyclic and so $H$ is nilpotent. This completes the proof.
\end{proof}
We conclude this section with the following interesting information on the group $M$ in the proof of Theorem \ref{th2}.
In fact, for the largest nilpotent group $M=\langle a,b\rangle$ relative to $a\in R_4(M)$, we have that $M/T$ is isomorphic to the largest (nilpotent) $2$-generated $4$-Engel group $E(2,4)$, where $T$ is the torsion subgroup of $M$ which is a $\{2,3,5\}$-group. Therefore, in a nilpotent $\{2,3,5\}'$-group, a right $4$-Engel element with an arbitrary element generate a $4$-Engel group. This can be seen by  comparing the presentations of $M/T$ and $E(2,4)$ as follows. One can obtain two finitely presented groups {\sf G1} and {\sf G2} isomorphic to $M/T$ and $E(2,4)$, respectively by {\sf GAP}:
\begin{verbatim}
MoverT:=FactorGroup(M,TorsionSubgroup(M));
E24:=NilpotentEngelQuotient(FreeGroup(2),4);
iso1:=IsomorphismFpGroup(MoverT);iso2:=IsomorphismFpGroup(E24);
G1:=Image(iso1);G2:=Image(iso2);
\end{verbatim}
Next, we find  the relators  of the groups {\sf G1} and {\sf G2} which are two sets of relators on 13 generators by the following command in {\sf GAP}.
\begin{verbatim}
r1:=RelatorsOfFpGroup(G1);r2:=RelatorsOfFpGroup(G2);
\end{verbatim}
Now, save these two sets of relators  by {\sf LogTo} command of {\sf GAP} in a file and go to the file to delete the terms as
\begin{verbatim}
<identity ...>
\end{verbatim}
in the sets {\sf r1} and {\sf r2}. Now call these two modified sets {\sf R1} and {\sf R2}. We show that {\sf R1=R2} as two sets of elements of the  free group {\sf f}  on 13 generators {\sf f1,f2,...,f13}.
\begin{verbatim}
f:=FreeGroup(13);
f1:=f.1;f2:=f.2;f3:=f.3;f4:=f.4;f5:=f.5;f6:=f.6;
f7:=f.7;f8:=f.8;f9:=f.9;f10:=f.11;f12:=f.12;f13:=f.13;
\end{verbatim}
Now by {\sf Read} function, load the file in {\sf GAP} and type the  simple command
{\sf R1=R2}. This gives us {\sf true} which shows $G_1$ and $G_2$ are two finitely presented groups with the same relators and generators and so they are isomorphic. We do not know if there is a guarantee that if someone else does as we did, then he/she finds the same relators for {\sf Fp} groups {\sf G1} and {\sf G2}, as we have found. Also we remark that using function {\sf IsomorphismGroups} to test if $G_1\cong G_2$, did not give us a result in less than 10 hours and we do not know whether this function can give us a result or not. \\

We summarize the above discussion as following.
\begin{thm}
Let $G$ be a nilpotent group generated by two elements, one of which is  a right $4$-Engel element. If $G$ has no element of order $2$, $3$ or $5$, then $G$ is a $4$-Engel group of class at most $6$.
\end{thm}
\section{\bf Right $n$-Engel elements for $n\geq 5$}\label{se4}
In this section we show that for every $n\geq 5$ there is a
nilpotent group $G$ of class $n+2$ containing elements $a$
and $x\in R_n(G)$  such that both $[x^{k},_n a]$ and $[x^{-1}, _na]$  have infinite order for all integers $k\geq 2$.

 Note that by Nickel's example \cite{nick2}, for every $n\geq 3$ we have already had a nilpotent group $K$ of class $n+2$ containing a right $n$-Engel element $x$ such that $[x^{-1},_n y]=[x^{2},_n y]\not=1$  for some $y\in K$ i.e, neither $x^2$ nor $x^{-1}$ are right $n$-Engel.  We have checked by {\sf nq} package of Nickel in {\sf GAP} that  $[x^{-1},_n y]=[x^{2},_n y]$ is of finite order whenever $n\in\{5,6,7,8\}$. In fact,
\begin{enumerate}
\item $o([x^{-1},_5 y])=3$, ~~~~~~~~  NqRuntime=1.7 Sec
\item $o([x^{-1},_6 y])=7$, ~~~~~~~~  NqRuntime=54.8 Sec
\item $o([x^{-1},_7 y])=4$, ~~~~~~~~  NqRuntime=1702 Sec
\item $o([x^{-1},_8 y])=9$, ~~~~~~~~ NqRuntime=56406 Sec
\end{enumerate}

Newman and Nickel \cite{newman} constructed a group $H$ as follows.
Let $F$ be the relatively free group, generated by $\{a,b\}$ with
nilpotency class $n+2$ and $\gamma_4( F)$ abelian. Let  $M$ be the
(normal) subgroup of $F$ generated by all commutators  in $a$, $b$
with at least  3 entries $b$ and the commutators $[b,_{n+1}a]$ and
$[b,_na,b]$. Then $\displaystyle H=\frac{F}{M}$. Note that the
normal closure of $b$ in $H$  is nilpotent of class 2.

We denote the generators of $H$ by $a, b$ again. Put $$ t=[b,_n a],
~~~u_j=[b,_{n-1-j}a,b,_j a],~~~ 0\leq j\leq n-2,$$ $$~~
u=\prod_{j=0}^{n-2}u_j,~~~v=[u_{n-2},a],~~~
w=\prod_{j=0}^{n-3}[u_j,a]
$$ and let $N$ be the subgroup $\langle tuw, t^2w, uw\rangle$.
Then $aN$ is a right $n$-Engel element in $\displaystyle \frac{H}{N}$ and
$[b,_n a]N$ has infinite order in $\displaystyle \frac{H}{N}$.

Now let $H$ be the above group and $N_0:=\langle
u,vw,vt^{-1}\rangle$. First, note that $N_0$ is a normal subgroup of $H$. For, clearly $t,v,w\in Z(H)$ and $u^b=u$. Also it is not hard to see that
$u_j^a=u_j[u_j,a]$ and thus $u^a=u vw$. This means that $N_0^a=N_0$ and so $N_0$ is a normal subgroup of $H$. Now we can state our main result of this section:
\begin{thm}\label{th3}
$[b,_n a]N_0=[b^{-2},_n a]N_0$ and it has infinite order in $\displaystyle
\frac{H}{N_0}$ and $[b^{-1},_n h]\in N_0$ for all $h\in H$.
Furthermore  $[b^{-k},_n a]N_0=v^{\binom{k}{2}}N_0$ for all $k\geq 2$.
\end{thm}
\begin{rem} \label{re1} As in \cite{newman}, the proof of  Theorem \ref{th3} involves a
series of commutator calculations based, as
usual, on the basic identities as following, which are mentioned in \cite{newman}. We bring them here for reader's convenience.
\begin{enumerate}
    \item $[g, cd]=[g,d][g, c] [g, c,d]$.
    \item $[cd,g]=[c,g][c,g,d][d,g]$.
    \item $[c^{-1},d]=[c,d,c^{-1}]^{-1}[c,d]^{-1}$.
    \item $[c,d^{-1}]=[c,d,d^{-1}]^{-1}[c,d]^{-1}$.
    \item $[hk,h_1,\dots,h_s]=[h,h_1,\ldots,h_s]$ for every $k$ in
    $\gamma_{n+3-s}(H)$ and arbitrary $h_1,\dots,h_s\in H$
    \item $[g,d,c]=[g,c,d][g,[d,c]]k$, where $k$ is a
    product of commutators of weight at least $4$ with entries $g$,
    $c$ and $d$.
    \item $[a,_n hk]=[a,_n h]$ for all  $h\in H$ and $k\in \gamma_3(H)$.
    \item $[g,d^{\delta}]=[g,d]^{\delta}[g,_2 d]^{(_2^{\delta})}k$,
    where $k$ is a product of commutators with at least $3$
    entries $d$ and $\delta$ is positive.
\end{enumerate}
\end{rem}
\noindent{\bf Proof of Theorem \ref{th3}.} By Remark \ref{re1}(7),
we may assume that $h$ is of the form
$a^{\alpha}b^{\beta}[b,a]^{\gamma}$. The following calculations
may depend to  the signs of $\alpha$ and $\beta$; we here
outline  only the case in which $\alpha$ and $\beta$ are positive.
\begin{eqnarray}
[b^{-1},_n h]&=&[b^{-1},_n a^{\alpha}b^{\beta}[b,a]^{\gamma}]\nonumber\\
&=&\displaystyle [b^{-1},_n a^{\alpha}b^{\beta}]
\prod_{j=0}^{n-1}[b^{-1},_{n-1-j}a^{\alpha}b^{\beta},[b,a]^{\gamma},_j a^{\alpha}b^{\beta}]\nonumber\\
&=&\displaystyle [b^{-1},_n a^{\alpha}b^{\beta}]\big([b,[b,a],_{n-1}a]
\prod_{j=0}^{n-2}[b,_{n-1-j}a,[b,a],_j
a]\big)^{-\alpha^{n-1}\gamma}.\nonumber
\end{eqnarray}
Since
\begin{eqnarray}
[b,[b,a],_{n-1}a]&=&[[[b,a],b]^{-1},_{n-1} a]\nonumber\\
&=&[b,a,b,_{n-1} a]^{-1}\nonumber\\
&=&v^{-1}\nonumber
\end{eqnarray}
and by Remark \ref{re1} (5) and (6)
$$[b,_{n-1-j}a,[b,a],_j a]=[b,_{n-j}a,b,_j
a]^{-1}[b,_{n-1-j}a,b,_{j+1}a]$$ we have
\begin{eqnarray}
\displaystyle\prod_{j=0}^{n-2}[b,_{n-1-j}a,[b,a],_j
a]&=&\prod_{j=0}^{n-2}[b,_{n-j}a,b,_j a]^{-1}[b,_{n-1-j}a,b,_{j+1} a]\nonumber\\
&=&\prod_{j=0}^{n-3}[b,_{n-1-j}a,b,_{j+1}
a]^{-1}\prod_{j=0}^{n-2}[b,_{n-1-j}a,b,_{j+1} a]
\nonumber\\
&=&v.\nonumber
\end{eqnarray}
Therefore
\begin{eqnarray}
\displaystyle[b^{-1},_n a^{\alpha}b^{\beta}[b,a]^{\gamma}]
&=&[b^{-1},_n a^{\alpha}b^{\beta}](v^{-1}v)^{-\alpha^{n-1}\gamma}\nonumber\\
&=&[b^{-1},_n a^{\alpha}b^{\beta}].\nonumber
\end{eqnarray}
On the other hand by Remark \ref{re1} (8) we have
\begin{eqnarray}
\displaystyle[b^{-1},_n a^{\alpha}b^{\beta}]&=& [b^{-1},_n a^{\alpha}]
\prod_{j=0}^{n-2}[b^{-1},_{n-1-j} a^{\alpha},b^{\beta},_j a^{\alpha}]\nonumber\\
&=&[b^{-1},_n a]^{\alpha^n}[b^{-1},_{n+1}
a]^{n\big{(}_2^{\alpha}\big{)}\alpha^{n-1}}(\displaystyle\prod_{j=0}^{n-2}
[b,_{n-1-j}a,b,_ja])^{-\alpha^{n-1}\beta}\nonumber\\
&&\times (\prod_{j=0}^{n-3}[b,_{n-1-j}
a,b,_{j+1}a])^{-(n-2)\big{(}_2^{\alpha}\big{)}\alpha^{n-2}\beta}\nonumber\\
&&\times[b,a,b,_{n-2}a]^{-(n-2)\big{(}_2^{\alpha}\big{)}\alpha^{n-2}\beta}
[b,_na,b]^{-(n-2)\big{(}_2^{\alpha}\big{)}\alpha^{n-2}\beta}\nonumber\\
&=&(vt^{-1})^{\alpha^n}
u^{-\alpha^{n-1}\beta}(vw)^{-(n-2)\big{(}_2^{\alpha}\big{)}\alpha^{n-2}\beta}.\nonumber
\end{eqnarray}
Therefore $b^{-1}N_0$ is a right $n$-Engel element in $\displaystyle
\frac{H}{N_0}$. This completes the second part of the theorem.

Since $\langle t,u,v,w\rangle$ is a free abelian group of rank 4,
it is clear that $[b,_na]N_0$ has infinite order. On the other hand
\begin{eqnarray}
\displaystyle [b^{-2},_na]&=&[[b^{-1}, a][b^{-1},a,b^{-1}][b^{-1},a],_{n-1}a]\nonumber\\
&=&[b^{-1},_n
a][b^{-1},a,b^{-1},_{n-1}a][b^{-1},a,b^{-1},[b,a],_{n-2}a][b^{-1},_na]\nonumber\\
&\equiv&[b,a,b,_{n-1}a]\mod N_0\nonumber\\
 &\equiv&v\mod N_0.\nonumber
\end{eqnarray}
Since $vt^{-1}\in N_0$ we have $[b,_n a]N_0=tN_0=vN_0=[b^{-2},_n a]N_0$. Now
let $k\geq 2$, $f(1)=0$ and $f(k)=(k-1)+f(k-1)=\binom{k}{2}$. Then
\begin{eqnarray}
\displaystyle [b^{-k},_na]&=&[[b^{-1}, a][b^{-1},a,b^{-(k-1)}][b^{-1},a],_{n-1}a]\nonumber\\
&=&[b^{-1},_n a][b^{-1},a,b^{-(k-1)},_{n-1}a][b^{-1},a,b^{-(k-1)},[b,a],_{n-2}a]
[b^{-(k-1)},_na]\nonumber\\
&\equiv&[b,a,b^{(k-1)},_{n-1}a]v^{f(k-1)}\mod N_0\nonumber\\
&\equiv&v^{f(k)}\mod N_0.\nonumber
\end{eqnarray}
This completes the proof. $\hfill \Box$\\

Now we answer negatively Question \ref{q1} which has been proposed in  \cite{ab1}.

Let $T$ be the torsion subgroup of $H/N_0$ and $x=bN_0T$ and $y=aN_0T$.  Then the group $\mathcal{M}=H/N_0T=\langle x,y\rangle$ is a torsion free, nilpotent of class $n+2$, $x\in R_n(\mathcal{M})$ and both $[x^{-1},_n y]$ and $[x^{k},_n a]$ are of infinite order for all integers $k\geq 2$. Since, for any given prime number $p$,  a finitely generated torsion-free nilpotent group is residually finite $p$-group, it follows that for any prime number $p$ and integer $k\geq 2$, there is a finite $p$-group $G(p,k)$ of class $n+2$ containing a right $n$-Engel element $t$ such that both $t^{k}$ and $t^{-1}$ are not  right $n$-Engel. This answers negatively Question \ref{q1}.

\section{\bf Subgroupness of the set of (bounded) Left Engel elements of a group}
Let $n=2^k\geq 2^{48}$ and $B(X,n)$ be the free Burnside  group on the set $X=\{x_i \;|\; i\in \mathbb{N}\}$ of the Burnside variety of exponent $n$ defined by the law
$x^n=1$. Lemma 6 of \cite{IO} states that  the subgroup $\langle x_{2k-1}^{n/2}x_{2k}^{n/2}\;|\; k=1,2,\dots\rangle$ of $B(X,n)$ is isomorphic  to $B(X,n)$ under the map $x_{2k-1}^{n/2}x_{2k}^{n/2}\rightarrow x_k$, $k=1,2,\dots$. Therefore  the subgroup
$\mathcal{G}:=\langle x_1^{n/2},x_2^{n/2},x_3^{n/2},x_4^{n/2}\rangle$ is generated by four elements of order $2$, contains the subgroup $\mathcal{H}=\langle x_1^{n/2}x_2^{n/2},x_3^{n/2}x_4^{n/2} \rangle$ isomorphic to the free $2$-generator Burnside group $B(2,n)$ of exponent $n$. One knows the tricky formulae $$[x,_ky]=[x,y]^{(-1)^{k-1}2^{k-1}}$$ holding for all elements $x$ and all elements $y$ of order $2$ in any group and all integers $k\geq 1$. It follows that the group $\mathcal{G}$ can be generated by four left $49$-Engel elements of $\mathcal{G}$. Thus $$\mathcal{G}=\langle L_{49}\big(\langle \mathcal{G}\rangle\big)\rangle=\langle L\big(\langle \mathcal{G}\rangle\big)\rangle=\langle \overline{L}\big(\langle \mathcal{G}\rangle\big)\rangle,$$
where $L(H)$ ($\overline{L}(H)$, resp.) denotes the set of (bounded, resp.) left Engel elements of a group $H$.\\

Suppose, if possible, $\mathcal{G}$ is an Engel group. Then $\mathcal{H}$ is also an Engel group. Let $Z$ and $Y$ be two free generators of $\mathcal{H}$. Thus $[Z,_k Y]=1$ for some integer $k\geq 1$. Since $\mathcal{H}$ is the free 2-generator Burnside group of exponent $n$, we have that every group of exponent $n$ is a $k$-Engel group. Therefore, $\mathcal{G}$ is an infinite finitely generated $k$-Engel group of exponent $n$, as $\mathcal{H}$ is infinite by a celebrated result of Ivanov \cite{I}. Hence, we have proved that
 \begin{prop}
 At least one of the following  happens.
 \begin{enumerate}
 \item There is an infinite finitely generated $k$-Engel group of exponent $n$ for some positive integer $k$ and $2$-power number $n$.
 \item There is a group $G$ such that  $L(G)=\overline{L}(G)$ and $L(G)$ is not a subgroup of $G$.
 \end{enumerate}
 \end{prop}
We believe that the subgroup $\mathcal{H}$ cannot be an Engel group, but we are unable to prove it.

\end{document}